\documentclass[11pt]{article}
\usepackage{amsfonts}
\title{AN OPTIMAL DIFFERENTIABLE SPHERE THEOREM FOR COMPLETE
MANIFOLDS \footnote{
\indent Research supported by the NSFC, Grant No.
10771187; the Trans-Century Training Programme
\newline \indent Foundation for Talents by the Ministry of
Education of China.}}
\author{HONG-WEI XU AND JUAN-RU GU}
\date{}
\begin{document}
\maketitle
\begin{abstract}
A new differentiable sphere theorem is obtained from the view of
submanifold geometry. An important scalar is defined by the scalar
curvature and the mean curvature of an oriented complete submanifold
$M^n$ in a space form $F^{n+p}(c)$ with $c\ge0$. Making use of the
Hamilton-Brendle-Schoen convergence result for Ricci flow and the
Lawson-Simons-Xin formula for the nonexistence of stable currents,
we prove that if the infimum of this scalar is positive, then $M$ is
diffeomorphic to $S^n$. We then introduce an intrinsic invariant
$I(M)$ for oriented complete Riemannian $n$-manifold $M$ via the
scalar, and prove that if $I(M)>0$, then $M$ is diffeomorphic to
$S^n$. It should be emphasized that our differentiable sphere
theorem is optimal for arbitrary $n(\ge2)$.
\end{abstract}
\section*{1. Introduction}
\hspace*{5mm}The investigation of curvature and topology of
Riemannian manifolds or submanifolds is one of the main stream in
global differential geometry. In 1898, Hadamard \cite{Hadamard}
proved a classical sphere theorem which says that any oriented
compact surface with positive Gaussian curvature in $R^3$ must be
diffeomorphic a sphere. It was seen from the Gauss-Bonnet theorem
that Hadamard's sphere theorem could be extended to the cases of
compact Riemannian surfaces with positive curvature. A natural
problem is stated as follows.\\\\
\textbf{Problem 1.1.} \emph{Is it possible to generalize the
Hadamard sphere theorem for compact Riemannian
surfaces to higher dimensional cases?}\\\\
\hspace*{5mm} In 1951, Rauch \cite{Rauch} first proved a topological
sphere theorem for positive pinched compact manifolds. During the
past six decades, there are many important progresses on topological
and differentiable pinching problems for Riemannian manifolds. The
most famous topological sphere theorem is Berger-Klingenberg's
quarter pinching theorem, which has been improved by many geometers
\cite{Berger,Brendle1,Shiohama2,Yau}. Recently Brendle and Schoen
\cite{Brendle3} obtained a classification of compact and simply
connected manifolds with weakly $1/4$-pinched curvatures.
Consequently, they obtained the following striking result.\\\\
\textbf{Theorem A.} \emph{Let M be an n-dimensional complete and
simply connected Riemannian manifold such that $1/4\le K_M\le1$.
Then $M$ is either diffeomorphic to $S^n$, or isometric to a compact
rank one symmetric
space (CROSS).}\\\\
\hspace*{5mm}Since the dimension of a complex projective space is
always even, Brendle and Schoen's differentiable sphere theorem is
optimal for even dimensional cases. More recently Petersen and Tao
\cite{Petersen1} have improved Brendle and Schoen's pinching
constant in Theorem A to $\frac{1}{4}-\varepsilon_n$, where
$\varepsilon_n$ is a positive constant depending only on $n$.
However, Petersen and Tao's pinching constant is not yet optimal for
odd dimensional
cases. \\
\hspace*{5mm}Let $M^{n}$ be an $n(\geq2)$-dimensional submanifold in
an $(n+p)$-dimensional Riemannian manifolds $N^{n+p}$. Denote by $H$
and $S$ the mean curvature and the squared length of the second
fundamental form of $M$ respectively. Using nonexistence for stable
currents on compact submanifolds of a sphere and the generalized
Poincare conjecture for dimension $n(\geq5)$ proved by Smale, Lawson
and Simons \cite{Lawson2} proved that if $M^{n}(n\ge 5)$ is an
oriented compact submanifold in $S^{n+p}$, and if $S<2\sqrt{n-1}$, then $M$ is homeomorphic to a sphere.\\
\hspace*{5mm}Let $F^{n+p}(c)$ be an $(n+p)$-dimensional simply
connected space form with nonnegative constant curvature $c$. Putting
$$\alpha(n,H,c)=nc+ \frac{n^{3}}{2(n-1)}H^{2} -
\frac{n(n-2)}{2(n-1)}\sqrt{n^{2}H^{4}+4(n-1)cH^{2}},$$ we have
$\min_{H}\alpha(n,H,c)=2\sqrt{n-1}c.$ Motivated by a rigidity
theorem in \cite{Xu90, Xu},
Shiohama and Xu \cite{Shiohama} improved Lawson-Simons' result and proved the following\\\\
\textbf{Theorem B.} \emph{Let $M^{n}(n\ge4)$ be an oriented complete
submanifold in $F^{n+p}(c)$ with $c\geq 0$. Suppose that
$\sup_{M}(S-\alpha(n,H,c))<0.$ Then $M$ is homeomorphic
to a sphere.}\\\\
\hspace*{5mm}The following
differentiable sphere theorem for hypersurfaces follows from the
convergence results for the mean curvature
flow and parabolic flow due to Huisken \cite{Huisken1} and Andrews \cite{Andrews}, respectively.\\\\
\textbf{Theorem C.} \emph{Let $M^{n}$ be an $n$-dimensional oriented
closed hypersurface in $F^{n+1}(c)$ with $c>0$. If $
S<2c+\frac{n^2H^2}{n-1},$ then $M$ is diffeomorphic to $S^n$.}\\\\
\hspace*{5mm}Recently, Xu and Zhao \cite{Xu2} proved a
differentiable sphere
theorem for submanifolds of a sphere with codimension $p(\ge1)$.\\\\
\textbf{Theorem D.} \emph{Let $M^n$ be an $n(\geq4)$-dimensional
oriented complete submanifold in $F^{n+p}(c)$ with $c>0$. Then \\
$(i)$ if $4\le n\le6$ and $\sup_{M}S<2\sqrt{n-1}c$, then $M$ is
diffeomorphic to $S^n$, \\
$(ii)$ if $n\geq7$ and $S<2\sqrt{2}c$, then $M$ is diffeomorphic to
$S^n$.}\\\\
\hspace*{5mm}Motivated by Theorems B, C and D, we propose the
following differentiable pinching problem.\\\\
\textbf{Problem 1.2.} \emph{Let $M^{n}$ be an oriented complete
submanifold in $F^{n+p}(c)$ with $c\ge0$. Suppose that $\sup_{M}
\Big(S-\frac{n^2H^{2}}{n-1}-2c\Big)<0$. Is it possible to prove that
$M$ must be diffeomorphic
to $S^{n}$?}\\\\
\hspace*{5mm}The purpose of the present article is to solve Problems
1.1 and 1.2, and prove some new differentiable pinching theorems for
complete submanifolds and Riemannian manifolds via Ricci flow and
stable currents. More precisely, we obtain the
following\\\\
\textbf{Theorem 1.1.} \emph{Let $M^n$ be an $n$-dimensional complete
submanifold in an $(n + p)$-dimensional Riemannian manifold
$N^{n+p}$.  Denote by $\overline{K}(\pi)$ the the sectional
curvature of N for tangent 2-plane $\pi(\subset T_xN)$ at point
$x\in N$. Set $\overline{K}_{\max}:=\max_{\pi\subset
T_{x}N}\overline{K}(\pi)$, $\overline{K}_{\min}:=\min_{\pi\subset
T_{x}N}\overline{K}(\pi)$. If
$S<\frac{8}{3}\Big(\overline{K}_{\min}-\frac{1}{4}\overline{K}_{\max}\Big)+
\frac{n^2H^{2}}{n-1}$,
 then $M$ is diffeomorphic to a space form. In particular, if M is simply connected, then M is diffeomorphic
to $S^n$ or $R^n$.}\\\\
\textbf{Theorem 1.2.} \emph{Let $M^n$ be an $n$-dimensional oriented
complete submanifold in $F^{n+p}(c)$ with $c\ge0$. If
$$\lambda(M):=\sup_M\Big(S- \frac{n^2H^{2}}{n-1}-2c\Big)<0,$$
 then $M$ is diffeomorphic to $S^n$.}\\\\
\hspace*{5mm} We shall show in Example 4.1 that Theorem 1.2 is
optimal for arbitrary $n(\ge2)$ and $p$. It follows from the Gauss
equation that the pinching condition in Theorem 1.2 is equivalent to
$$\mu(M):=-\lambda(M)=\inf_M\Big[R-
\frac{n^2(n-2)}{n-1}H^{2}-(n+1)(n-2)c\Big]>0,$$ where $R$ is the
scalar curvature of $M$.\\
\hspace*{5mm} For a complete Riemannian $n$-manifold $M^n$,  we set
$\mathcal{C}:=\{\varphi \, ;\varphi:M\longrightarrow F^{n+p}(c)$ is
an isometric embedding for some constant $c\ge0$ and positive
integer $p$ $\}$ and $\mathcal{D}:=\{\varphi \,
;\varphi:M\longrightarrow R^{n+p}$ is an isometric embedding for
some positive integer $p$ $\}.$ With the aid of the Nash embedding
theorem, we get $\mathcal{C}\supset \mathcal{D}\neq\emptyset$. We
define two intrinsic invariants $I(M)$ and $I_0(M)$ by
$$I(M):=\sup_{\varphi\in\mathcal{C}}\mu(M, \varphi):=\sup_{\varphi\in\mathcal{C}}
\inf_M\Big[R- \frac{n^2(n-2)}{n-1}H^{2}-(n+1)(n-2)c\Big],$$
$$I_0(M):=\sup_{\varphi\in\mathcal{D}}\mu(M, \varphi):=\sup_{\varphi\in\mathcal{D}}
\inf_M\Big[R- \frac{n^2(n-2)}{n-1}H^{2}\Big].$$
Notice that $I(M)\ge I_0(M)$. We shall prove\\\\
\textbf{Theorem 1.3.} \emph{Let $M^n$ be an oriented complete
Riemannian $n$-manifold. If $I(M)>0,$ then $M$ is diffeomorphic to
$S^n$. In particular, if $I_0(M)>0,$ then $M$ is
diffeomorphic to $S^n$.}\\\\
\textbf{Remark 1.1.} In the case $n=2$, Theorem 1.3 is reduce to the
Hadamard sphere theorem for compact Riemannian surfaces. We shall
give an example (Example 4.2) to show that our differentiable sphere
theorem for Riemannian manifolds is optimal
for arbitrary $n(\ge2)$.\\\\
\hspace*{5mm}Furthermore, we obtain some other differentiable
pinching theorems for complete submanifolds in Riemannian manifolds,
which extend the sphere theorems due to Huisken, Xu and Zhao
\cite{Huisken1, Xu2}.\\
\section*{2. Notation and lemmas}
\hspace*{5mm}Let $M^{n}$ be an $n(\geq2)$-dimensional submanifold in
an $(n+p)$-dimensional Riemannian manifolds $N^{n+p}$. We shall make
use of the following convention on the range of indices.
$$ 1\leq A,B,C,\ldots\leq n+p;\ 1\leq i,j,k,\ldots\leq n;\ n+1\leq
\alpha,\beta,\gamma,\ldots\leq n+p.$$ For an arbitrary fixed point
$x\in M\subset N$, we choose an orthonormal local frame field
$\{e_{A}\}$ in $N^{n+p}$ such that $e_{i}$'s are tangent to $M$.
Denote by $\{\omega_{A}\}$ the dual frame field of $\{e_{A}\}$. Let
$Rm$ and $\overline{Rm}$ be the Riemannian curvature tensors of $M$
and $N$ respectively, and $h$ the second fundamental form of $M$.
Then

$$Rm=\sum_{i,j,k,l}{R_{ijkl}}\omega_{i}\otimes\omega_{j}\otimes\omega_{k}\otimes\omega_{l},\hspace*{13mm}$$
$$\overline{Rm}=\sum_{A,B,C,D}\overline{R}_{ABCD}\omega_{A}\otimes\omega_{B}\otimes\omega_{C}\otimes\omega_{D},$$
$$h=\sum_{\alpha,i,j}h^{\alpha}_{ij}\omega_{i}\otimes\omega_{j}\otimes
e_{\alpha},\hspace*{29mm}$$
\begin{equation}
R_{ijkl}=\overline{R}_{ijkl}+\sum_{\alpha}(h^{\alpha}_{ik}h^{\alpha}_{jl}-h^{\alpha}_{il}h^{\alpha}_{jk}).\hspace*{11mm}
\end{equation} The squared norm $S$ of the second fundamental form
and the mean curvature $H$ of $M$ are given by
$$
S:=\sum_{\alpha,i,j}(h^{\alpha}_{ij})^{2},\,\,H:=\Big|\frac{1}{n}\sum_{\alpha,i}h^{\alpha}_{ii}e_{\alpha}\Big|.$$
Denote by $K(\pi)$ the the sectional curvature of $M$ for tangent
2-plane $\pi(\subset T_xM)$ at point $x\in M$, $\overline{K}(\pi)$
the the sectional curvature of $N$ for tangent 2-plane $\pi(\subset
T_xN)$ at point $x\in N$. Set $\overline{K}_{\min}:=\min_{\pi\subset
T_{x}N}\overline{K}(\pi)$, $\overline{K}_{\max}:=\max_{\pi\subset
T_{x}N}\overline{K}(\pi)$.\\
\hspace*{5mm}The Lawson-Simons-Xin non-existence theorem
\cite{Lawson2, Xin} for stable currents in a compact Riemannian
manifold $M$ isometrically immersed into $F^{n+p}(c)$ is employed to
eliminate the homology groups
$H_{q}(M;Z)$ for $0<q<n$.\\\\
\textbf{Lemma 2.1.} \emph{Let $M^{n}$ be a compact submanifold in
$F^{n+p}(c)$ with $c\geq 0$. Assume that$$
\sum_{k=q+1}^{n}\sum_{i=1}^{q}[2|h(e_{i},e_{k})|^{2}-\langle
h(e_{i},e_{i}),h(e_{k},e_{k})\rangle]<q(n-q)c$$ holds for any
orthonormal basis $\{e_{i}\}$ of $M_{x}$ at any point $x\in M$,
where q is an integer satisfying $0<q<n$. Then there does not exist
any stable q-currents. Moreover,
$$H_{q}(M;Z)=H_{n-q}(M;Z)=0,$$ where $H_{i}(M;Z)$ is the $i$-th homology group of M with integer
coefficients.}\\\\
 \hspace*{5mm}The following convergence result for the Ricci flow, initialed by Brendle and Schoen \cite{Brendle2}, was finally obtained by Brendle \cite{Brendle}.\\\\
\textbf{Lemma 2.2.}\emph{ Let $(M,g_0)$ be a compact Riemannian
manifold of dimension $n(\geq4)$. Assume that
\begin{equation}R_{1313}+\lambda^2 R_{1414} + R_{2323} + \lambda^2R_{2424}-
2\lambda R_{1234} > 0\end{equation} for all orthonormal four-frames
$\{e_1,e_2,e_3,e_4\}$ and all $\lambda\in[-1,1]$. Then the
normalized Ricci flow with initial metric $g_0$
$$\frac{\partial}{\partial t}g(t) = -2Ric_{g(t)} +\frac{2}{n}
r_{g(t)}g(t),$$ exists for all time and converges to a constant
curvature metric as $t\rightarrow\infty$. Here $r_{g(t)}$ denotes
the mean value of the scalar curvature of $g(t)$.}\\
\section*{3. Sphere Theorem in dimension three}
\hspace*{5mm}Using Lemma 2.1 and the assumption for $S$, we
obtain the following\\\\
\textbf{Theorem 3.1.} \emph{Let $M^{3}$ be a $3$-dimensional
oriented compact submanifold in a simply connected space form
$F^{3+p}(c)$ with nonnegative constant curvature c.
 If $S<2c+\frac{9}{2}H^2$, then $M$ is
diffeomorphic to $S^3$.}\\\\
\textbf{Proof.} We observe that
\begin{eqnarray}
&&\sum_{k=q+1}^{3}\sum_{i=1}^{q}[2|h(e_{i},e_{k})|^{2}-\langle
h(e_{i},e_{i}),h(e_{k},e_{k})\rangle]\nonumber \\
&=&2\sum_{\alpha}\sum_{k=q+1}^{3}\sum_{i=1}^{q}(h^{\alpha}_{ik})^{2}-
\sum_{\alpha}\sum_{k=q+1}^{3}\sum_{i=1}^{q}h^{\alpha}_{ii}h^{\alpha}_{kk} \nonumber \\
&=&\sum_{\alpha}\Big[2\sum_{k=q+1}^{3}\sum_{i=1}^{q}(h^{\alpha}_{ik})^{2}-
\Big(\sum_{i=1}^{q}h^{\alpha}_{ii}\Big)\Big(\sum_{i=1}^{3}h^{\alpha}_{ii}-\sum_{i=1}^{q}h^{\alpha}_{ii}\Big)\Big].
\end{eqnarray}
Setting
$$S_{\alpha}:=\sum_{i,j=1}^{3}(h^{\alpha}_{ij})^{2},\ \
 T_{\alpha}:=\sum_{i=1}^{3}h^{\alpha}_{ii} ,\ \
 \tilde{S}_{\alpha}:=\sum_{i=1}^{3}(h^{\alpha}_{ii})^{2},$$ we have
$$S=\sum_{\alpha}S_{\alpha},\ \
 9H^{2}=\sum_{\alpha}T_{\alpha}^{2},$$
 and
 \begin{equation}qr\tilde{S}_{\alpha}=qr\sum_{i=1}^{q}(h^{\alpha}_{ii})^{2}+qr\sum_{k=q+1}^{3}(h^{\alpha}_{kk})^{2}
\geq
 r\Big(\sum_{i=1}^{q}h^{\alpha}_{ii}\Big)^{2}+q\Big(\sum_{k=q+1}^{3}h^{\alpha}_{kk}\Big)^{2},\end{equation}
where $r:=3-q$. Inserting
$$T_{\alpha}-\sum_{i=1}^{q}h^{\alpha}_{ii}=\sum_{k=q+1}^{3}h^{\alpha}_{kk},$$
into the right hand side of (4), we
get\begin{equation}3\Big(\sum_{i=1}^{q}h^{\alpha}_{ii}\Big)^{2}-
2qT_{\alpha}\sum_{i=1}^{q}h^{\alpha}_{ii}+qT_{\alpha}^{2}-qr\tilde{S}_{\alpha}\leq
0.\end{equation} Set
$$Z_{\alpha}:=-\Big(\sum_{i=1}^{q}h^{\alpha}_{ii}\Big)\Big(T_{\alpha}-\sum_{i=1}^{q}h^{\alpha}_{ii}\Big).$$
It follows from (5) that
\begin{equation}3Z_{\alpha}+(r-q)T_{\alpha}\sum_{i=1}^{q}h^{\alpha}_{ii}+qT^{2}_{\alpha}-qr\tilde{S}_{\alpha}\leq 0.\end{equation}
Making use of the relations
$$\sum_{i=1}^{3}\Big(h^{\alpha}_{ii}-\frac{T_{\alpha}}{3}\Big)^{2}=\tilde{S}_{\alpha}-\frac{T_{\alpha}^{2}}{3},\
\ \sum_{i=1}^{3}\Big(h^{\alpha}_{ii}-\frac{T_{\alpha}}{3}\Big)=0,\ \
\sum_{i=1}^{q}\Big(h^{\alpha}_{ii}-\frac{T_{\alpha}}{3}\Big)+\frac{q}{3}T_{\alpha}=
\sum_{i=1}^{q}h^{\alpha}_{ii},$$ and setting
$\tilde{h}^{\alpha}_{ii}:=h^{\alpha}_{ii}-\frac{T_{\alpha}}{3}$, we
obtain
\begin{equation}
\tilde{S}_{\alpha}-\frac{T_{\alpha}^{2}}{3}\geq\frac{1}{q}\Big(\sum_{i=1}^{q}\tilde{h}^{\alpha}_{ii}\Big)^{2}
+\frac{1}{r}\Big(\sum_{k=q+1}^{3}\tilde{h}^{\alpha}_{kk}\Big)^{2}
=\Big(\frac{1}{q}+\frac{1}{r}\Big)\Big[\sum_{i=1}^{q}\Big(h^{\alpha}_{ii}-\frac{T_{\alpha}}{3}\Big)\Big]^2.\end{equation}
Therefore we find
\begin{equation}\Big|\sum_{i=1}^{q}\Big(h^{\alpha}_{ii}-\frac{T_{\alpha}}{3}\Big)\Big|
\leq\sqrt{\frac{qr}{3}\Big(\tilde{S}_{\alpha}-\frac{T_{\alpha}^{2}}{3}\Big)}.\end{equation}
This together with (6) implies
\begin{equation}Z_{\alpha}\leq\frac{qr}{3}\tilde{S}_{\alpha}-\Big[\frac{q(r-q)}{9}
+\frac{q}{3}\Big]T_{\alpha}^{2}+\frac{|r-q|}{3}|T_{\alpha}|\sqrt{\frac{qr}{3}
\Big(\tilde{S}_{\alpha}-\frac{T_{\alpha}^{2}}{3}\Big)}.\end{equation} From
(3), (9) and the fact $qr=2$ and $|r-q|=1$, we obtain
\begin{eqnarray}
&&\sum_{k=q+1}^{3}\sum_{i=1}^{q}[2|h(e_{i},e_{k})|^{2}-\langle
h(e_{i},e_{i}),h(e_{k},e_{k})\rangle]-qrc \nonumber \\
&\leq&\sum_{\alpha}\Big[S_{\alpha}-\frac{\tilde{S}_{\alpha}}{3}-\frac{4}{9}T_{\alpha}^{2}+
\frac{|T_{\alpha}|}{3}\sqrt{\frac{2}{3}\Big(\tilde{S}_{\alpha}-\frac{T_{\alpha}^{2}}{3}\Big)}\Big]-2c \nonumber \\
&\leq&S-4H^{2}-2c-\frac{1}{3}\sum_{\alpha}\tilde{S}_{\alpha}+\sum_{\alpha}\Big[\frac{T_{\alpha}^{2}}{18}+\frac{1}{3}\Big(\tilde{S}_{\alpha}-\frac{T_{\alpha}^{2}}{3}\Big)\Big] \nonumber \\
&=&S-\frac{9}{2}H^{2}-2c
\end{eqnarray}
Then under the assumption, we obtain
\begin{equation}
\sum_{k=q+1}^{3}\sum_{i=1}^{q}[2|h(e_{i},e_{k})|^{2}-\langle
h(e_{i},e_{i}),h(e_{k},e_{k})\rangle] -qrc<0.\end{equation} Suppose
that $\pi_{1}(M)\neq0$. Since $M$ is compact, it follows from a
classical theorem due to Cartan and Hadamard that there exists a
minimal closed geodesic in any non-trivial homotopy class in
$\pi_{1}(M)$. However, combining Lemma 2.1 and (11) we know that
there does not exist any stable integral currents on $M$. This
contradicts the hypothesis. Therefore,
$\pi_{1}(M)=0$. \\
\hspace*{5mm} It follows from Proposition 2.1 in \cite{Shiohama} and
the assumption for $S$ that
\begin{eqnarray*}Ric_{M}(X)&\ge&\frac{2}{3}\Big[3c+6H^2-S-\frac{3}{\sqrt{6}}H(S-3H^2)^{1/2}\Big]\\
&=&\frac{2}{3}\Big[(3c+\frac{27}{4}H^2-\frac{3}{2}S)+\frac{3}{4}H^2+\frac{1}{2}(S-3H^2)
-\frac{3}{\sqrt{6}}H(S-3H^2)^{1/2}\Big]>0\end{eqnarray*}
holds for any unit vector $X\in T_x{M}.$ By Hamilton's convergence
result for Ricci flow in three dimensions \cite{Hamilton}, it
follows that $M$ is diffeomorphic to a 3-dimensional spherical
space form. This completes the proof of Theorem 3.1.\\\\
\textbf{Corollary 3.1.}\emph{ Let $M$ be a $3$-dimensional oriented
compact submanifold in the unit sphere $S^{3+p}$. Suppose that $H\ge
\frac{2}{3}\sqrt{\sqrt2-1}$.  If $S<2\sqrt{2}$,
 then $M$ is diffeomorphic to
$S^3$.}\\\\
\textbf{Proof.} By a direct computation, we get
$$S<2+\frac{9}{2}H^2.$$
By Theorem 3.1, we see that $M$ is
diffeomorphic to $S^3$. This proves Corollary 3.1.\\\\
\hspace*{5mm}Up to now, the following problem proposed by Lawson
and Simons \cite{Lawson2} is still open.\\\\
\textbf{Problem 3.1.} \emph{Let $M$ be a $3$-dimensional oriented
compact submanifold in the unit sphere $S^{3+p}$. Suppose that
$S<2\sqrt{2}$. Can one prove that $M$ must be diffeomorphic to
$S^3$?}\\
\section*{4. Differentiable sphere theorem in higher dimensions}
\hspace*{5mm}The following lemma will be used in the proof of our theorems.\\\\
\textbf{Lemma 4.1.} \emph{Let $M^{n}$ be an $n$-dimensional
 submanifold in an
$(n+p)$-dimensional Riemannian manifold $N^{n+p}$, and $\pi$ a
tangent 2-plane on $T_xM$ at point $x\in M$. Choose an orthonormal
two-frame $\{e_1,e_2\}$ at $x$ such that $\pi=span\{e_1,e_2\}$. Then
\begin{equation}K(\pi)\geq\frac{1}{2}\Big(2\overline{K}_{\min}+\frac{n^2H^2}{n-1}-S\Big)+\sum_{\alpha=n+1}^{n+p}\sum_{j>i,(i,j)\neq(1,2)}(h^{\alpha}_{ij})^2.
\end{equation}}\\
\textbf{Proof.} We extend the orthonormal two-frame $\{e_1,e_2\}$ to
$\{e_{1},\cdots,e_{n+p}\}$ such that $e_{i}$'s are tangent to $M$.
Setting $S_{\alpha}:= \sum_{i,j=1}^{n}(h^{\alpha}_{ij})^{2},$ we
have
\begin{equation}
 \Big(\sum_{i=1}^nh_{ii}^{\alpha}\Big)^2=(n-1)\Big[\sum_{i=1}^{n}(h^{\alpha}_{ii})^2+\sum_{i\neq
j}(h_{ij}^{\alpha})^2+\frac{(\sum_{i=1}^nh_{ii}^{\alpha})^2}{n-1}-S_{\alpha}\Big].\end{equation}
Note that
\begin{eqnarray*}
\Big(\sum_{i=1}^nh_{ii}^{\alpha}\Big)^2&\leq&(n-1)\Big[(h^{\alpha}_{11}+h^{\alpha}_{22})^2+\sum_{i>2}(h^{\alpha}_{ii})^2\Big]\\
&=&(n-1)\Big[\sum_{i=1}^{n}(h^{\alpha}_{ii})^2+2h^{\alpha}_{11}h^{\alpha}_{22}\Big].\end{eqnarray*}
This together with (13) implies
\begin{equation}
2h^{\alpha}_{11}h^{\alpha}_{22}\geq\sum_{i\neq
j}(h_{ij}^{\alpha})^2+\frac{(\sum_{i=1}^nh_{ii}^{\alpha})^2}{n-1}-S_{\alpha}.\end{equation}
From the Gauss equation and (14) we get\begin{eqnarray}
K(\pi)&=&\overline{R}_{1212}+\sum_{\alpha=n+1}^{n+p}[h^{\alpha}_{11}h^{\alpha}_{22}-(h^{\alpha}_{12})^2] \nonumber \\
&\geq&\sum_{\alpha=n+1}^{n+p}\Big[\sum_{j>2}(h_{1j}^{\alpha})^2+\sum_{j>2}(h_{2j}^{\alpha})^2+\frac{1}{2}\sum_{i\neq
j>2}(h_{ij}^{\alpha})^2\Big]+\frac{1}{2}\Big(\frac{n^2H^2}{n-1}-S\Big)+\overline{K}_{\min} \nonumber \\
 &\geq&\frac{1}{2}\Big(2\overline{K}_{\min}+\frac{n^2H^2}{n-1}-S\Big)+\sum_{\alpha=n+1}^{n+p}\sum_{
  j>i,(i,j)\neq(1,2)}(h^{\alpha}_{ij})^2. \end{eqnarray}
This proves Lemma 4.1.\\\\
\textbf{Lemma 4.2.} \emph{Let $M^{n}$ be an $n$-dimensional complete
submanifold in an $(n+p)$-dimensional Riemannian manifold $N^{n+p}$.
If $\sup_M\Big(S-2\overline{K}_{\min}-\frac{n^2H^2}{n-1}\Big)<0$,
then
M is compact.}\\\\
\textbf{Proof.} From the assumption and Lemma 4.1, it follows that
there exists an $\varepsilon>0$ such that $K_{M}\geq\varepsilon$.
 By the Bonnet-Myers's theorem, we know that $M$ is compact. This completes the proof.\\\\
\textbf{Theorem 4.1.} \emph{Let $(M,g_0)$ be an
$n(\geq4)$-dimensional complete submanifold in an $(n +
p)$-dimensional Riemannian manifold $N^{n+p}$. If
$\sup_M\Big[S-\frac{8}{3}\Big(\overline{K}_{\min}-\frac{1}{4}\overline{K}_{\max}\Big)-
\frac{n^2H^{2}}{n-1}\Big]<0$,
 then the normalized Ricci flow with initial metric $g_0$
$$\frac{\partial}{\partial t}g(t) = -2Ric_{g(t)} +\frac{2}{n}
r_{g(t)}g(t),$$ exists for all time and converges to a constant
curvature metric as $t\rightarrow\infty$. Moreover, $M$ is
diffeomorphic to a space form. In particular, if M is simply
connected, then M is diffeomorphic
to $S^n$.}\\\\
\textbf{Proof.} By Lemma 4.2, it follows that $M$ is compact. When
$n\geq4$, suppose $\{e_1,e_2,e_3,e_4\}$ is an orthonormal four-frame
and $\lambda\in R$. From the Gauss equation (1) and Berger's inequality
we have
\begin{eqnarray}|R_{1234}|&=&|\overline{R}_{1234}+\sum_\alpha(h^{\alpha}_{13}h^{\alpha}_{24}-h^{\alpha}_{14}h^{\alpha}_{23})|\nonumber \\
&\leq&\frac{2}{3}(\overline{K}_{\max}-\overline{K}_{\min})+\sum_\alpha|h^{\alpha}_{13}h^{\alpha}_{24}-h^{\alpha}_{14}h^{\alpha}_{23}|.\end{eqnarray}
This together with Lemma 4.1 implies
\begin{eqnarray}&&R_{1313}+\lambda^2 R_{1414} + R_{2323} +
\lambda^2R_{2424}- 2\lambda R_{1234} \nonumber \\
&\geq&(1+\lambda^2)\Big(2\overline{K}_{\min}+\frac{n^2H^2}{n-1}-S\Big) \nonumber \\
&&+\sum_{\alpha}\sum_{i<j,(i,j)\neq(1,3)}(h^{\alpha}_{ij})^2+\sum_{\alpha}\sum_{i<j,(i,j)\neq(2,3)}(h^{\alpha}_{ij})^2 \nonumber \\
&&+\lambda^2\Big[\sum_{\alpha}\sum_{i<
j,(i,j)\neq(1,4)}(h^{\alpha}_{ij})^2+\sum_{\alpha}\sum_{i<
j,(i,j)\neq(2,4)}(h^{\alpha}_{ij})^2\Big] \nonumber \\
&&-2|\lambda|\Big[\frac{2}{3}(\overline{K}_{\max}-\overline{K}_{\min})+\sum_\alpha|h^{\alpha}_{13}h^{\alpha}_{24}-h^{\alpha}_{14}h^{\alpha}_{23}|\Big] \nonumber \\
&\geq&(1+\lambda^2)\Big[\frac{8}{3}\Big(\overline{K}_{\min}-\frac{1}{4}\overline{K}_{\max}\Big)+\frac{n^2H^2}{n-1}-S\Big] \nonumber \\
&&+\sum_{\alpha}[(h_{24}^{\alpha})^2+\lambda^2(h_{13}^{\alpha})^2+(h_{14}^{\alpha})^2+\lambda^2(h_{23}^{\alpha})^2] \nonumber \\
&&-2|\lambda||h^{\alpha}_{13}h^{\alpha}_{24}|-2|\lambda||h^{\alpha}_{14}h^{\alpha}_{23}| \nonumber \\
&\geq&(1+\lambda^2)\Big[\frac{8}{3}\Big(\overline{K}_{\min}-\frac{1}{4}\overline{K}_{\max}\Big)+\frac{n^2H^2}{n-1}-S\Big] \nonumber \\
&>&0.\end{eqnarray} It follows from Lemma 2.2 that $M$ is
diffeomorphic to a space form. In particular, if $M$ is simply
connected, then $M$ is diffeomorphic
to $S^n$. This completes the proof Theorem 4.1.\\\\
\textbf{Proof of Theorem 1.1.} By the assumption, we have
\begin{eqnarray}S&<&\frac{8}{3}\Big(\overline{K}_{\min}-\frac{1}{4}\overline{K}_{\max}\Big)+
\frac{n^2H^{2}}{n-1} \nonumber \\
&\le&2\overline{K}_{\min}+\frac{n^2H^2}{n-1}.\end{eqnarray} So
$$2\overline{K}_{\min}+\frac{n^2H^2}{n-1}-S>0.$$
This together with Lemma 4.1 implies $K_{M}>0$. \\
\hspace*{5mm}When $M$ is non-compact, a theorem due to
Cheeger-Gromoll-Meyer \cite{Cheeger,Gromoll3} says
that $M$ must be diffeomorphic to $R^n$. \\
\hspace*{5mm}When $M$ is compact, we consider the following cases:
(i) If $n=2$, it follows from the fact $K_{M}>0$ that $M$ is
diffeomorphic to $S^2$ or $RP^2$. (iii) If $n=3$, Hamilton's theorem
\cite{Hamilton} shows that $M$ is diffeomorphic to a spherical space
form. (iii) If $n\geq 4$, the assertion follows from Theorem 4.1. In
particular, when $M$ is simply connected, we conclude that $M$ must
be diffeomorphic to $S^n$ or $R^n$. This completes
the proof of Theorem 1.1.\\\\
\textbf{Theorem 4.2.}\emph{ Let $M^n$ be an $n(\geq4)$-dimensional
oriented complete submanifold in an $(n + p)$-dimensional simply
connected space form $F^{n+p}(c)$ with nonnegative constant
curvature $c$. If $\sup_M\Big(S-\frac{n^2H^{2}}{n-1}\Big)<2c$,
then $M$ is diffeomorphic to $S^n$.}\\\\
\textbf{Proof.} It is easy to see that
\begin{eqnarray}\alpha(n,H,c)&=&nc+\frac{n^{3}}{2(n-1)}H^{2}
-\frac{n(n-2)}{2(n-1)}\sqrt{n^{2}H^{4}+4(n-1)cH^{2}} \nonumber \\
&\geq&nc+\frac{n^{3}}{2(n-1)}H^{2}
-\frac{n(n-2)}{2(n-1)}\Big[\frac{nH^2}{2}+\frac{n^{2}H^{2}+4(n-1)c}{2n}\Big] \nonumber \\
&=&2c+\frac{n^2H^2}{n-1}.
\end{eqnarray}
It follows from Theorem B that $M$ is a topological
sphere.\\
\hspace*{5mm}On the other hand, we see from Theorem 4.1 that $M$ is
diffeomorphic to a space form. Therefore, $M$ is
diffeomorphic to $S^n$. This proves Theorem 4.2.\\\\
\textbf{Theorem 4.3.}\emph{ Let $M^n$ be an $n$-dimensional oriented
complete submanifold in an $(n + p)$-dimensional simply connected
space form $F^{n+p}(c)$ with nonnegative constant curvature $c$. If
$S<2c+ \frac{n^2H^{2}}{n-1}$,
 then $M$ is diffeomorphic to $S^n$ or $R^n$.}\\\\
\textbf{Proof.} From Lemma 4.1, we know that $K_{M}>0$. When $M$ is non-compact, the assertion
follows from the proof of Theorem 1.1.\\
\hspace*{5mm}When $M$ is compact, we consider the following cases:
(i) If $n=2$, from the Gauss-Bonnet theorem we see that the genus of
$M$ is zero, and hence $M$ is a topological sphere. Therefore, $M$
is diffeomorphic to $S^2$. (ii) If $n\geq 3$, it follows from
Theorems 3.1 and 4.2 that $M$ is diffeomorphic to $S^n$.
This completes the proof.\\\\
\textbf{Proof of Theorem 1.2.} From Lemma 4.2, we know that $M$ is
compact. This together with Theorem 4.3 implies
that $M$ is diffeomorphic to $S^n$. This completes the proof of Theorem 1.2.\\\\
\hspace*{5mm}The following example shows that the pinching
conditions in Theorems
1.2 and 4.3 are the best possible for arbitrary $n(\ge2)$ and $p$.\\\\
\textbf{Example 4.1.} (i) When $c=0$, let
$M:=S^{n-1}\Big(\frac{n-1}{nH_0}\Big)\times R^1\subset
R^{n+1}\subset R^{n+p}$, where $H_0$ is a positive constant. Then
$H=H_0$ and $S=\frac{n^2H^2}{n-1}$. (ii) When $c>0,$ without loss of
generality, we only consider the case $c=1$. Let
$M:=S^{1}\Big(\frac{1}{\sqrt{1+\lambda^2}}\Big)\times
S^{n-1}\Big(\frac{\lambda}{\sqrt{1+\lambda^2}}\Big)\subset
S^{n+1}\subset S^{n+p}$, where $\lambda$ is a positive constant. We
have $H=\frac{1}{n}[\lambda-(n-1)\frac{1}{\lambda}]$ and
$S=\lambda^2+(n-1)\frac{1}{\lambda^2}$. Then
$S-\frac{n^2H^2}{n-1}-2=\frac{(n-2)}{(n-1)}\lambda^2.$ Thus, for any
$\varepsilon>0$ we can find a submanifold
$M:=S^{1}\Big(\frac{1}{\sqrt{1+\lambda^2}}\Big)\times
S^{n-1}\Big(\frac{\lambda}{\sqrt{1+\lambda^2}}\Big)\subset S^{n+p}$
satisfying $S<2+\frac{n^2H^2}{n-1}+\varepsilon.$\\\\
\textbf{Proof of Theorem 1.3.} By the assumption, we put
$I(M):=\varepsilon_{0}>0.$ There exists an isometric embedding
$\varphi:M\longrightarrow F^{n+p}(c)$ such that
$$\mu(M, \varphi)\ge \frac{1}{2}\varepsilon_{0}>0.$$ Thus $\lambda(M, \varphi)<0.$ It follows from Theorem 1.2 that $M$ is
diffeomorphic to  $S^n$. This proves Theorem 1.3.\\\\
\hspace*{5mm}The following example shows that Theorem 1.3 is optimal
for
arbitrary $n(\ge2)$.\\\\
\textbf{Example 4.2.} Let
$M:=S^{n-1}\Big(\frac{n-1}{nH_0}\Big)\times R^1\subset
R^{n+1}\subset R^{n+p}$, where $H_0$ is a positive constant. We
consider the inclusion $\varphi_0:M\longrightarrow R^{n+p}.$
Following Example 4.1, we have
$$\mu(M,
\varphi_0)=-\lambda(M, \varphi_0)=0.$$ This implies
 $I(M)\ge I_{0}(M)\ge0.$ Since $M$ is not diffeomorphic to $S^n$, it
follows from Theorem 1.3
that $I_{0}(M)\le I(M)\le0.$ Hence $I(M)=I_{0}(M)=0.$\\\\
\hspace*{5mm}Finally we present the following differentiable sphere
theorem for even dimensional submanifolds in a general Riemannian
manifold, which is an extension of Theorems 1.2 and 4.3 as well as
the sphere theorems due to Huisken, Xu and Zhao
\cite{Huisken1, Xu2}.\\\\
\textbf{Theorem 4.4.} \emph{Let $M^n$ be an even dimensional
oriented complete submanifold in an $(n + p)$-dimensional Riemannian
manifold $N^{n+p}$. Then \\
$(i)$ if
$S<\frac{8}{3}\Big(\overline{K}_{\min}-\frac{1}{4}\overline{K}_{\max}\Big)+
\frac{n^2H^{2}}{n-1}$, then M is diffeomorphic to $S^n$ or $R^n$.\\
$(ii)$ if
 $\sup_{M} \Big[S-\frac{8}{3}\Big(\overline{K}_{\min}-\frac{1}{4}\overline{K}_{\max}\Big)-
\frac{n^2H^{2}}{n-1}\Big]<0$, then M is diffeomorphic to $S^n$.}\\\\
\textbf{Proof.} (i) It follows from the assumption and Lemma 4.1
that $K_{M}>0$. When $M$ is non-compact, it follows from
Cheeger-Gromoll-Meyer's soul theorem \cite{Cheeger,Gromoll3} that
$M$ is diffeomorphic to $R^n$. When $M$ is compact, it's seen from
the assumption and Synge's theorem that $M$ is simply connected.
This together with Theorem 1.1 implies that $M$ is diffeomorphic to
$S^n$. Therefore, we conclude that $M$ is diffeomorphic
to $S^n$ or $R^n$.\\
\hspace*{5mm} (ii) From the assumption and Lemma 4.2, we see that
$M$ is compact. This together with (i) implies that $M$ is
diffeomorphic to $S^n$. This completes the proof of Theorem 4.4.\\\\
\textbf{\large{Acknowledgement.}}\ The authors would like to thank
Professors Kefeng Liu and Richard Schoen for their
helpful discussions and valuable suggestions.\\

Center of Mathematical Sciences\

Zhejiang University\

Hangzhou 310027\

China\\\

E-mail address: xuhw@cms.zju.edu.cn; gujuanru@sina.com

\end{document}